%
%
%
%
%
%
%
%
%
%
%
%
%
%
%

\magnification=\magstephalf
\input amstex
\documentstyle{amsppt}
\NoBlackBoxes

\def\sqr#1#2{{\vcenter{\hrule height.#2pt
        \hbox{\vrule width.#2pt height#1pt \kern#1pt
                \vrule width.#2pt}
        \hrule height.#2pt}}}
\def\square{\mathchoice\sqr64\sqr64\sqr{4}3\sqr{3}3}

\def\QED{\hfill$\square$}
\def\s{{\bold s}}
\def\p{{\frak p}}
\def\m{{\frak m}}

\topmatter

\title
A Generalized Dedekind--Mertens Lemma And \\ Its Converse
\endtitle

\rightheadtext{A Generalized Dedekind--Mertens Lemma And Its Converse}

\author
Alberto Corso,
William Heinzer, and
Craig Huneke
\endauthor

\address
Department of Mathematics, Purdue University,
West Lafayette, IN 47906 USA
\endaddress

\email
corso, heinzer, and huneke\@math.purdue.edu
\endemail

\date
\
\enddate

\dedicatory
To Wolmer V. Vasconcelos on the occasion of his sixtieth birthday
\enddedicatory

\thanks
The authors gratefully acknowledge partial support from the NSF.
\endthanks

\keywords
Content ideal of a polynomial, integral closure of an ideal,
Dedekind--Mertens Lemma, Dedekind--Mertens number
\endkeywords

\subjclass
Primary: 13A15; Secondary: 13B25, 13G05, 13H10
\endsubjclass

\abstract
We study content ideals of polynomials and their behavior under multiplication.
We give a generalization of the Lemma of Dedekind--Mertens and prove
the converse under suitable dimensionality restrictions.
\endabstract

\endtopmatter

\document

\head
1. Introduction
\endhead

\baselineskip 15 pt

Let $R$ be a commutative ring and let $t$ be an indeterminate over $R$. For a
polynomial $f \in R[t]$, the {\it  content ideal} $c(f)$ of $f$ is the
ideal of $R$ generated by the coefficients of $f$.

If $f$ and $g$ are two polynomials in $R[t]$ one clearly has that
$c(fg) \subseteq c(f)c(g)$ and the classical Lemma of Gauss, in one of its
forms, says that equality holds if $R$ is a principal ideal domain.
More generally, $c(f)c(g)$ is integral over
$c(fg)$. Since ideals in principal ideal domains are integrally closed,
Gauss's Lemma follows from this statement. An even more precise statement is
given by the Lemma of Dedekind--Mertens.
The lemma asserts that if $f$ and $g$ are two polynomials and $n$ is
the degree of $g$ then
$$  
c(fg) c(f)^n = c(f)^{n+1}c(g). \tag 1
$$
Interchanging the roles of $f$ and $g$, there is obviously an analogous
formula involving the degree of the polynomial $f$ and powers of $c(g)$.

There has recently been renewed interest in this lemma for a variety of
reasons. See \cite{AG}, \cite{AK}, \cite{BG}, \cite{CVV},
\cite{GGP}, \cite{GV}, \cite{HH1}, \cite{HH2}, and \cite{N}.

Easy examples show that the degrees of the
polynomials may be a too crude measure of the relation between the ideals
$c(f)c(g)$ and $c(fg)$.
In order to obtain a sharper form of (1) as well as a finer measure
of comparison between $c(f)c(g)$ and $c(fg)$, the last two authors of
the present paper introduced the {\it Dedekind--Mertens number} of a
polynomial $g \in R[t]$ (see \cite{HH2}). This number $\mu_R(g)$ is defined
as the smallest positive integer $k$ such that
$$  
c(fg) c(f)^{k-1} = c(f)^k c(g) \tag 2
$$
for every polynomial $f\in R[t]$.

The relation between the minimal number of generators $\mu_R(c(g))$ of
$c(g)$ and $\mu_R(g)$ is addressed in \cite{HH2},  and the main result of
that paper states that
$$
\mu_R(g) \leq \mu_R(c(g)) \tag 3
$$
(see \cite{HH2, Theorem 2.1}). Since $\mu_R(c(g))\leq \text{deg}(g) + 1$,
this statement implies the usual Dedekind--Mertens Lemma.
In the same paper, the following question is raised (see
\cite{HH2, Question 1.3}): Let $(R, \m)$ be an excellent local domain, and let
$g \in R[t]$. Is $\mu(g) = \mu(c(g))$?
In the case $\mu(g)=1$, the veracity of the above equality reduces to
a question posed in the early sixties in the Ph.D. thesis of
I. Kaplansky's student H.T. Tsang. An affirmative answer to the
question of Tsang in a broad variety of cases (including all
Noetherian domains) is given in two recent works on Gaussian
polynomials: see \cite{GV} and \cite{HH1}.

The main theorem in this paper proves the converse under one extra assumption
on the dimension of the ring. Theorem~4.2 includes all of
Theorem~1.1 and a bit more.

\proclaim{Theorem~1.1}  
Let $(R, \m)$ be a universally catenary, analytically unramified
Noetherian local ring. Suppose $g \in R[t]$ has Dedekind--Mertens number
$\mu_R(g) = k$. Assume that $\dim (R/\p) \geq k$ for all minimal
primes $\p$ of $R$. Then $\mu(c(g)) \leq k$. Therefore, $\mu_R(c(g)) =
\mu_R(g)$.
\endproclaim

In Section 5 we  prove that at least in certain examples, the
assumption of Theorem~1.1 concerning the dimension is necessary.
We give in Example~5.1 a polynomial over a one-dimensional complete
local Gorenstein domain whose content ideal requires three generators,
but whose Dedekind--Mertens number is two.

In order to prove Theorem~1.1, we introduce another
version of the Dedekind--Mertens number which we call the
polarized Dedekind--Mertens number.
We define the {\it polarized Dedekind--Mertens number} $\widetilde{\mu}_R(g)$
of a polynomial $g \in R[t]$ to be the smallest positive integer $k$ such that
$$  
\sum_{i=1}^k c(f_ig) c(f_1) \cdots \widehat{c(f_i)} \cdots c(f_k)
=
c(f_1) \cdots c(f_k) c(g) \tag 4
$$
for all polynomials $f_1, \ldots, f_k \in R[t]$, where $\widehat{c(f_i)}$
indicates the deletion of $c(f_i)$.
We  refer to (4) as the {\it polarized Dedekind--Mertens
formula}. Clearly, (2) follows from (4) by choosing
$f = f_1 = \ldots = f_k$; thus we have:

\remark{Remark~1.2}  
For every polynomial $g \in R[t]$, the Dedekind--Mertens number
$\mu_R(g)$ is less than or equal to $\widetilde{\mu}_R(g)$, the polarized
Dedekind--Mertens number of $g$.
\endremark

\medskip

It turns out that we rely heavily on the {\it a priori } stronger version
of Dedekind--Mertens provided by this polarized form. We
show in Theorem~2.5 that the polarized Dedekind--Mertens number
$\widetilde{\mu}_R(g)$ is related to the minimal number of generators of $c(g)$
in the same way as the Dedekind--Mertens number $\mu_R(g)$, namely
$$
\widetilde{\mu}_R(g) \leq \mu_R(c(g)).
$$
This raises the issue of whether $\mu_R(g) = \widetilde{\mu}_R(g)$.
In Theorem~2.8 we establish this equality under certain dimensionality
restrictions.

The strategy for the proof of our main theorem is as follows:
we first study the $0$-dimensional local Gorenstein case in
detail. To study the polarized form of the Dedekind--Mertens lemma in such
a ring, it is convenient to study its dual form, which effectively converts
information about equality of ideals into linear equations. In
particular the following theorem is an important step in this
translation:

\proclaim{Theorem~3.1}
Let $(R,\m)$ be a local Artinian ring, let $\s = (0 \colon \m)$ denote the
socle of $R$, and let $s$ be the dimension of $\s$ as a vector space over
$R/\m$.
Suppose $g(t) \in R[t]$ is a polynomial of degree $n$ and let
$I = (0 :_R c(g))$. For any ideal $J$ such that $I \subseteq J \subseteq
(I \colon \m)$ and $\dim(J/I)=r >s$ and for any $m$ such that
$(m+1)r > (n+m+1)s$ there exists a polynomial $f(t)$ of degree $m$ with the
following properties
\roster
\item"({\it a})"
$c(f) \subseteq J$;

\item"({\it b})"
$fg=0$;

\item"({\it c})"
$c(f)c(g)\not=0$.
\endroster
\endproclaim

\medskip

After our study of $0$-dimensional Gorenstein rings, we are ready
to give the proof of the main result. We reduce to the
$0$-dimensional Gorenstein case by using the fact that the rings
we are dealing with are {\it approximately Gorenstein}. This means
that one of the following two equivalent conditions holds (see \cite{Ho}):

\roster
\item"(i)"
For every integer $n>0$ there is an ideal $I \subseteq \m^n$ such that
$R/I$ is Gorenstein.

\item"(ii)"
For every integer $n>0$ there is an $\m$-primary irreducible ideal
$I \subseteq \m^n$.
\endroster

The statement of Theorem~4.2 is the same as the one of Theorem~1.1
with the additional conclusion that under these hypotheses we have
$$
\mu_R(c(g)) = \widetilde{\mu}_R(g) = \mu_R(g).
$$

The restriction on the dimension of $R/\p$ for minimal primes $\p$ of
$R$ comes into play from the translation into linear equations which
occurs in the $0$-dimensional case. At this point, we need to ensure
that we have more variables than equations, and after winding back
to our original situation, we need to know there exist ideals $I_1, \ldots,
I_k$ such that the number of minimal generators of the product of
these ideals is sufficiently larger than the minimal number of
generators of $\sum_i I_1\cdots \widehat{I_i} \cdots I_k$. To prove
this, we need to assume the dimension is sufficiently large. For example,
over a $1$-dimensional local ring, the number of generators of
an arbitrary ideal is bounded by the multiplicity of the ring.

In Section 5  several classes of examples are developed
over one-dimensional domains which show that
additional  assumptions are necessary in general to have the formula
$\mu_R(g) = \mu(c(g))$.

\head
2. Bounding the polarized Dedekind--Mertens number
\endhead

Throughout the paper we use the integral closure of an ideal.
Recall the definition:

\definition{Definition 2.1}
Given an ideal $I$ of a ring $R$, an element $x \in R$ is in the
{\it integral closure} $\overline I$ of
$I$  if $x$ satisfies an equation of the
form $x^k + a_1x^{k-1} + \cdots + a_k = 0$ where $a_i \in I^i$.
\enddefinition

\remark{Remark 2.2}
As mentioned above, it is well known that $c(f)c(g)\subseteq
\overline{c(fg)}$, so the ideals  $c(f)c(g)$ and $c(fg)$ have
the same integral closure. For example, see \cite{E}.
\endremark

\medskip

Another observation we use is:

\remark{Remark~2.3} 
For every ideal $I$ of a ring $R$ and polynomial $g \in R[t]$ the
Dedekind--Mertens number and polarized Dedekind--Mertens number of the
image of $g$ in $(R/I)[t]$ are less than or equal to the corresponding numbers
associated to $g$.
\endremark

\medskip

Let $R$ be a commutative ring and let $S$ be a polynomial extension
of $R$. A polynomial $g$ over $R$ can also be viewed as a polynomial
over $S$. Clearly, one has that $\mu_R(g) \leq \mu_S(g)$. Lemma~2.4
relates $\widetilde{\mu}_R(g)$ to $\mu_S(g)$, where $S$ is a suitable
polynomial extension of $R$.

\proclaim{Lemma~2.4}  
Let $R$ be a commutative ring.
If $g \in R[t]$ is such that for all polynomials $F \in S[t]$, with
$S=R[x_1, \ldots, x_{k-1}]$ and $x_1, \ldots, x_{k-1}$ indeterminates,
the equality
$$ 
c(Fg)c(F)^{k-1}=c(F)^k c(g) \tag 5
$$
holds, then $\widetilde{\mu}_R(g) \leq k$. In particular, if
$\mu_S(g) = k$, then
$(5)$ holds, so $\widetilde{\mu}_R(g) \leq \mu_S(g)$.
\endproclaim
\demo{Proof}
For any $f_1, \ldots, f_k \in R[t]$ consider the polynomial
$$
F=f_1+ \sum_{i=2}^k f_ix_{i-1}t^{N_{i-1}} \in R[x_1, \ldots, x_{k-1}][t]
$$
where the $N_i$'s are chosen recursively so that
$N_1 > \max \{ \deg_t(f_1), \deg_t(f_1g) \}$ and
$$
N_i > \max \{ \deg_t(f_i), \deg_t(f_ig) \} + N_{i-1}, \quad \text{\rm for \ }
i>1.
$$

By assumption the equality
$$
\multline
c\biggl(\biggl(f_1 + \sum_{i=2}^k f_ix_{i-1}t^{N_{i-1}}\biggr)g\biggr)\,
\biggl(c\biggl(f_1 + \sum_{i=2}^k
f_ix_{i-1}t^{N_{i-1}}\biggr)\biggr)^{k-1} = \\
\biggl(c\biggl(f_1 + \sum_{i=2}^k f_ix_{i-1}t^{N_{i-1}}\biggr)\biggr)^k
\, c(g)
\endmultline
$$
holds, and from our choice of the $N_i$'s it follows that this is equivalent to
$$
\multline
\biggl(c(f_1g) + \sum_{i=2}^k c(f_ig) x_{i-1} \biggr)\, \biggl(c(f_1) +
\sum_{i=2}^k c(f_i) x_{i-1} \biggr)^{k-1} = \\
= \biggl(c(f_1) + \sum_{i=2}^k c(f_i) x_{i-1} \biggr)^k \, c(g).
\endmultline
$$
Finally, the comparison of the coefficient of the monomial
$x_1 \cdots x_{k-1}$ in both terms of the previous identity gives the
formula
$$
\sum_{i=1}^k c(f_ig) c(f_1) \cdots \widehat{c(f_i)} \cdots c(f_k)
=
c(f_1) \cdots c(f_k) c(g)
$$
and thus shows that $\widetilde{\mu}_R(g) \leq k$.
\QED
\enddemo

\medskip

It is shown in \cite{HH2,Theorem 2.1} that the Dedekind--Mertens
number $\mu(g)$ is bounded above by the number of generators needed
to generate locally the content ideal of $g$. Applying this result
to a polynomial extension ring of $R$ yields the following
consequence of Lemma~2.4.

\proclaim{Theorem~2.5}  
For a polynomial  $g \in R[t]$, the polarized Dedekind--Mertens number
$\widetilde{\mu}_R(g)$ is bounded above by the number of local generators
of the content ideal $c(g)$ of $g$, i.e., we have
$\widetilde{\mu}_R(g) \le \mu_R(c(g))$.
\endproclaim

\remark{Remark~2.6}
Suppose $R$ is a subring of a ring $S$.  If
$S$ is faithfully flat over $R$ and
$g \in R[t] \subseteq S[t]$ is a polynomial, it is easily seen that
$\mu_R(g) \leq \mu_S(g)$. Sometimes, however,  this inequality
is strict. For example, if $(R,\m)$ is a local Artinian ring that is
not Gorenstein, then as noted in \cite{HH1, Remark 1.6},
there exists a polynomial $f \in R[t]$ such
that $\mu_R(f) = 1$, while $\mu(c(f)) > 1$. It
follows \cite{HH1, Remark 1.7} that if $S$ is a polynomial ring extension
in two variables over $R$, then  $\mu_S(f) > 1$.
\endremark

\remark{Remark~2.7}  
We use in several places the following result of Rees \cite{Re2, Theorem 2.1}:
Suppose $(R, \m)$ is a formally equidimensional local ring of dimension
$d$ and $I=(a_1, \ldots, a_r)$ is an ideal of the principal class, i.e.,
$\dim(R/I)=d-r$. If $F(X_1, \ldots, X_r) \in R[X_1, \ldots, X_r]$ is a
homogeneous polynomial such that  $F(a_1, \ldots, a_r)$ is nilpotent, then the
coefficients of $F$ belong to the integral closure $\overline{I}$ of $I$.
\endremark

\medskip

We show in Theorem~2.8 the equality of the Dedekind--Mertens number
and the polarized Dedekind--Mertens number under certain
dimensionality restrictions.

\proclaim{Theorem~2.8}  
Let $(R, \m)$ be a universally catenary, analytically unramified
Noetherian local ring.  Suppose $g \in R[t]$ has Dedekind--Mertens number
$\mu_R(g) = k$. If $\dim(R/\p) \ge k$ for each minimal prime $\p$ of $R$,
then $\widetilde{\mu}_R(g) \le k$.  Therefore
$\widetilde{\mu}_R(g) = \mu_R(g)$.
\endproclaim
\demo{Proof}
Let $f_1, \ldots, f_k$ be arbitrary polynomials in $R[t]$. For
simplicity, we
let $I = c(f_1) \cdots c(f_k)c(g)$ and
$J = \sum_{i=1}^k c(f_ig) c(f_1) \cdots \widehat{c(f_i)}
\cdots c(f_k)$.

By the Artin--Rees lemma, there  exists an integer $t$ so that
$\m^t \cap I \subseteq \m I$. Also, by \cite{Re1, Theorem 1.4} there
exists an integer $n$ with the property that $\overline{\m^n} \subseteq \m^t$,
where $\overline{\m^n}$ denotes the integral closure of $\m^n$. Choose now
$a_1, \ldots, a_k \in \m^n$ so that for all minimal primes $\p$ of $R$
the height of $((a_1, \ldots, a_k)R+\p)/\p$ is exactly $k$.

Consider the polynomial
$$
f = f_1a_1 + \sum_{i=2}^k f_ia_it^{N_{i-1}} \in R[t]
$$
where the $N_i$'s are chosen recursively so that
$N_1 > \max \{ \deg_t(f_1), \deg_t(f_1g) \}$ and
$$
N_i > \max \{ \deg_t(f_i), \deg_t(f_ig) \} + N_{i-1}, \quad \text{\rm for \ }
i>1.
$$

By assumption we have $c(fg)c(f)^{k-1} = c(f)^kc(g)$, so the equality
$$
\split
c\biggl(\biggl(f_1a_1+ \sum_{i=2}^k f_ia_it^{N_{i-1}}\biggr)g\biggr)\,
\biggl(c\biggl(
                 & f_1a_1+ \sum_{i=2}^k
                   f_ia_it^{N_{i-1}}\biggr)\biggr)^{k-1} = \\
                 & = \biggl(c\biggl(f_1a_1 + \sum_{i=2}^k
                   f_ia_it^{N_{i-1}}\biggr)\biggr)^k \, c(g)
\endsplit
$$
holds. From our choice of the $N_i$'s, this is equivalent to
$$
\split
\biggl(c(f_1g)a_1 + \sum_{i=2}^k c(f_ig) a_i \biggr)\, \biggl(c(f_1)a_1
+
                 & \sum_{i=2}^k c(f_i) a_i \biggr)^{k-1} = \\
                 & \biggl(c(f_1)a_1 + \sum_{i=2}^k c(f_i) a_i
                   \biggr)^k \, c(g).
\endsplit
$$
This last equality implies that if $K$ is the ideal of $R$
generated by the elements of the form $a_1^{e_1} \cdots  a_k^{e_k}$,
where the $e_i$ are nonnegative integers such that
$\sum_{i=1}^k e_i = k$ and not all $e_i = 1$, then
$$ 
I a_1 \cdots a_k
\subseteq
(K, J a_1 \cdots a_k). \tag 6
$$
Equation $(6)$ implies that for every $b \in I$ there
exists $c \in J$  such that
$$ 
a_1 \cdots a_k (b-c) \in K. \tag 7
$$
Since $R$ is universally catenary, $R/\p$ is formally equidimensional
for each minimal prime $\p$ of $R$ \cite{Mat, Theorem 31.7}.
Equation $(7)$ implies the existence of a homogeneous
polynomial $F(X_1, \dots , X_k) \in (R/\p)[X_1,\dots, X_k]$ of
degree $k$ having the image of
$(b-c)$ in $R/\p$ as the coefficient of $X_1\cdots X_k$ and having
the property that
$F(a_1,\dots ,a_k) = 0$.
Hence by the result of Rees in Remark~2.7, the image of
$b-c$ in $R/\p$ is integral over the image of $(a_1,\ldots , a_k)$
in $R/\p$ for every minimal prime $\p$. Thus
$$
b-c \in \overline{(a_1, \ldots, a_k)} \cap I \subseteq \overline{\m^n}
\cap I \subseteq \m^t \cap I \subseteq \m I,
$$
so that $I \subseteq J + \m I$, hence by Nakayama's lemma $I = J$.
\QED
\enddemo

\head
3. Zero-Dimensional results
\endhead

\proclaim{Theorem~3.1} 
Let $(R,\m)$ be a local Artinian ring, let $\s = (0 \colon \m)$ denote the
socle of $R$, and let $s$ be the dimension of $\s$ as a vector space over
$R/\m$.
Suppose $g(t) \in R[t]$ is a polynomial of degree $n$ and let
$I = (0 :_R c(g))$. For any ideal $J$ such that $I \subseteq J \subseteq
(I \colon \m)$ and $\dim(J/I)=r >s$ and for any $m$ such that
$(m+1)r > (n+m+1)s$ there exists a polynomial $f(t)$ of degree $m$ with the
following properties
\roster
\item"({\it a})"
$c(f) \subseteq J$;

\item"({\it b})"
$fg=0$;

\item"({\it c})"
$c(f)c(g)\not=0$.
\endroster
\endproclaim
\demo{Proof}
Let $g(t) = a_nt^n + a_{n-1}t^{n-1} + \hdots + a_1t + a_0 \in R[t]$.
The condition  for another polynomial
$f(t) = b_mt^m + \hdots + b_1t + b_0$ that the product $fg = 0$
is clearly equivalent to a system of $n+m+1$ linear equations in
the $m+1$ `variables' $b_m, \dots ,b_1,b_0$. More precisely,
for each $d$ with $0 \le d \le m+n$, we have
$$ 
\sum_{i+j = d}a_ib_j = 0. \tag 8
$$

Let $A$ denote the matrix of coefficients
of the system of $n+m+1$ equations in  $m+1$ variables
defined by $(8)$.
Let $x_1,\dots x_r$ be the preimages in $J$ of a basis for $J/I$ over $R/\m$,
and let $y_1,\dots , y_s$ be a basis of $\s$. As $a_ix_j \m =0$ for any
$0 \leq i \leq n$ and $1 \leq j \leq r$, there exist $c_{ijk} \in R$
such that
$$ 
a_ix_j = \sum_{k=1}^s c_{ijk}y_k . \tag 9
$$
For indeterminates
$\{z_{ij}\}_{ 0 \le i \le m,\, 1 \le j \le r}$ over $R$,
let $h(t) = w_mt^m + \dots + w_1t + w_0$, where
$$ 
w_i = z_{i1}x_1 + \dots + z_{ir}x_r
\qquad \text{ for } \quad 0 \le i \le m. \tag 10
$$
The product $gh = 0$ is equivalent to the system of linear equations:
$A \vec w = \vec 0$, where $\vec w$ is the transpose of $(w_m,\dots,w_0)$.

In general, for each $d$ with $0 \le d \le m+n$,
$\sum_{i+j = d}a_iw_j = 0$ implies, by $(10)$,
$$
\sum_{i+j = d}a_i(z_{j1}x_1 + \dots
+ z_{jr}x_r) = 0,
$$
which, by $(9)$, implies
$$ 
\left(\sum_{i+j = d}(c_{i1k}z_{j1} + \dots + c_{irk}z_{jr})\right)y_k  = 0,
\quad 0 \le d \le m+n, \  1\le k \le s. \tag 11
$$
Since $y_1, \dots, y_s$ form a basis for the socle $\s$  of $R$, to solve
the system of equations $(11)$ over $R$ in the $(m+1)r$ variables
$\{z_{ij}\}$ is equivalent to solving  over the residue field of $R$
the  system
$$ 
\sum_{i+j = d}(c_{i1k}'z_{j1} + \dots + c_{irk}'z_{jr})  = 0,
\qquad 0 \le d \le m+n, \quad 1\le k \le s,  \tag 12
$$
where $c_{ijk}'$ denotes the image in $R/\m$ of $c_{ijk} \in R$.

In $(12)$ we have a system of $(m+n+1)s$ homogeneous linear equations
over the field $R/\m$ in the $(m+1)r$ variables $\{z_{ij}\}$. For any $m$
such that $(m+1)r > (m+n+1)s$ the system $(12)$ has a nontrivial
solution $z_{ij} = e_{ij}'$. Let $e_{ij} \in R$ be a preimage of $e_{ij}'$
for $0 \le i \le m, 1 \le j \le r$.
Then for $f(t) = b_nt^n + \dots + b_1t + b_0$, where $b_i = e_{i1}x_1 +
\dots + e_{ir}x_r$, we have $fg = 0$ in $R[t]$. The linear independence of
the images of the $x_i$ in $J/I$ implies that $c(f)$ is not contained in $I$.
Hence $c(f)c(g) \ne 0$.
\QED
\enddemo

\medskip

Given $g(t) \in R[t]$ under what conditions do there exists  polynomials $f$
and $h$ in $R[t]$ such that $c(fg)c(h)=0=c(hg)c(f)$ but $c(f)c(g)c(h)\not=0$?
Theorem~3.3 below answers this question in a special case.

\remark{Remark~3.2}
Let $R$ be a zero-dimensional local Gorenstein ring.
An immediate consequence of duality is that for each  ideal  $I$ of $R$ the
dimension of the socle of $R/(0 \colon I)$ equals the minimal number of
generators of $I$ (see \cite{BH, Proposition 3.2.12}).
\endremark

\proclaim{Theorem~3.3} 
Suppose $(R, \m)$ is a zero-dimensional local Gorenstein ring and let
$g, f_1,$ $\ldots, f_n \in R[t]$. Set
$$
A = c(f_1) \cdots c(f_n) \qquad \text{and} \qquad
B = \sum_{i=1}^n c(f_ig) c(f_1) \cdots
\widehat{c(f_i)} \cdots c(f_n)
$$
where $\widehat{c(f_i)}$ indicates that $c(f_i)$ has been omitted.
If
$$ 
\dim(A c(g)/ (\m Ac(g) + B)) \geq \mu(A)+1, \tag 13
$$
then there exists an $h\in R[t]$ such that $A c(h) c(g) \not=0$ but
$B c(h) = 0 = c(hg) A$.
\endproclaim
\demo{Proof}
Let $J= 0 \colon (\m A c(g)+ B)$. Using duality
and the following inclusions
$$
0 \colon A c(g) \subseteq J \subseteq 0 \colon \m A c(g)
$$
one concludes that
$$
\dim(J/(0 \colon A c(g))) > \dim(\text{soc}(R/(0 \colon A))=\mu(A).
$$
By Theorem~3.1 applied to the ring $R/(0 \colon A)$
there exists $h(t)\in R[t]$ such that $c(h) \subseteq J$ and $A c(h)c(g)
\not= 0$ but $c(hg) A =0$.
\QED
\enddemo

\head
4. Minimal generators of the content ideal
\endhead

We observe in Proposition~4.1 good behavior of the
(polarized) Dedekind--Mertens number under passage from a
Noetherian local ring to its completion.

\proclaim{Proposition~4.1} 
Let $(R,\m)$ be a local Noetherian ring and let
$(\widehat{R}, \widehat{\m})$
denote the $\m$-adic completion of $R$.
For every  polynomial $g(t) \in R[t]$  we have
$$
\mu_R(g) = \mu_{\widehat{R}}(g) \quad \text{ and } \quad
\widetilde{\mu}_R(g) = \widetilde{\mu}_{\widehat{R}}(g).
$$
In other words, the $($polarized$)$ Dedekind--Mertens number
of $g$ over $R$ equals the $($polarized$)$ Dedekind--Mertens
number of $g$ over $\widehat{R}$.
\endproclaim
\demo{Proof}
Let $g^*$ denote the image of $g$ in $(R/\m^s)[t]$.
Since all ideals of $R$ are closed in the $\m$-adic topology on $R$,
$\mu_R(g) = \mu_{R/\m^s}(g^*)$  and
$\widetilde{\mu}_R(g) = \widetilde{\mu}_{R/\m^s}(g^*)$ for
sufficiently large $s$.
Moreover,  $R/\m^s \cong
\widehat{R}/\widehat{\m}^s$ for each $s$, so we have
$$
\mu_{R/\m^s}(g^*) = \mu_{\widehat{R}/\widehat{\m}^s}(g^*) \quad \text{ and }
\quad \widetilde{\mu}_{R/\m^s}(g^*) =
\widetilde{\mu}_{\widehat{R}/\widehat{\m}^s}(g^*),
$$
from which the assertions follow.
\QED
\enddemo

\proclaim{Theorem~4.2} 
Let $(R, \m)$ be a universally catenary,
analytically unramified Noetherian local ring.
Suppose $g \in R[t]$ has Dedekind--Mertens number
$\mu_R(g) = k$. Assume that $\dim (R/\p)
\geq k$ for all minimal primes $\p$ of $R$. Then
$\mu(c(g)) \leq k$. Therefore, $\mu_R(c(g)) =
\widetilde{\mu}_R(g)= \mu_R(g)$.
\endproclaim
\demo{Proof}
By Theorem~2.8,
the polarized Dedekind--Mertens number  $\widetilde{\mu}_R(g)=k$, i.e.,
for all polynomials $f_1, \ldots, f_k \in R[t]$  we have
$$
\sum_{i=1}^k c(f_ig) c(f_1) \cdots \widehat{c(f_i)} \cdots c(f_k)
=
c(f_1) \cdots c(f_k) c(g).
$$

In view of Proposition~4.1, we may assume that $R$ is
complete.
Suppose the ideal $c(g)$ is minimally generated by $z_1, \ldots, z_m$
with $m \geq k+1$.

By the Artin--Rees lemma, there exists an integer $t$ so that $c(g) \cap
\m^t \subseteq \m c(g)$. Also, by \cite{Re1, Theorem 1.4} there exists an
integer $n$ with the property that $\overline{\m^n} \subseteq \m^t$.

Choose now $a, b_1, \ldots b_{k-1} \in \m^n$ so that for all minimal primes
$\p$ of $R$ the height of $(a, b_1, \ldots, b_{k-1})R+\p/\p$ is exactly $k$
and consider the ideals
$$
\prod_{i=1}^{k-1} (a, b_i)^{r_i}, \qquad c(g)\prod_{i=1}^{k-1} (a, b_i)^{r_i},
$$
where the $r_i$'s are nonnegative integers.

For simplicity we use ${\bold b}^{\bold v}$
for $b_1^{v_1} \cdots b_{k-1}^{v_{k-1}}$, $|{\bold v}|$ to denote
$ \sum_{j=1}^{k-1} v_j$,  the length of ${\bold v}$, and $P$ to
denote  $\sum_{j=1}^{k-1} r_j$.

We  observe that the minimal number of generators of
$\prod_{i=1}^{k-1} (a, b_i)^{r_i}$ is
$$
\mu \left(\prod_{i=1}^{k-1} (a, b_i)^{r_i} \right) = \prod_{i=1}^{k-1} (r_i+1).
$$
Indeed, the ideal $\prod_{i=1}^{k-1} (a, b_i)^{r_i}$ is generated
by elements of
the form $a^u{\bold b}^{\bold v}$, where
$u+ |{\bold v}| = P$ and $0 \leq v_i \leq r_i$ for all $1 \leq i \leq
k-1$. For every minimal prime $\p$ of $R$,
Remark~2.7 implies the images of these elements in $R/\p$
form an irredundant generating set. Hence, these elements
are irredundant generators in $R$.
Moreover, these elements are in one-to-one correspondence with the $k$-tuples
$(u, v_1, \ldots, v_{k-1})$ of integers satisfying the above restrictions;
an easy calculation shows that the number of such $k$-tuples is exactly
$\prod_{i=1}^{k-1} (r_i+1)$.

We claim that the minimal number of generators of
$c(g)\prod_{i=1}^{k-1} (a, b_i)^{r_i}$ is exactly
$$
\mu \left( c(g)\prod_{i=1}^{k-1} (a, b_i)^{r_i} \right) =
m\prod_{i=1}^{k-1} (r_i+1).
$$
If not, there exists an element $z_ia^u {\bold b}^{\bold v}$,
with $u + |{\bold v}| = P$, which can be written as a combination of
the remaining ones, namely
$$ 
\multline
\qquad
z_ia^u {\bold b}^{\bold v} = z_i \Biggl( \sum \Sb u'+|{\bold v'}| = P
\\ {\bold v} \not = {\bold v'} \endSb  q_{(i,u',{\bold v'})}
a^{u'} {\bold b}^{\bold v'}\Biggr) + \\
\sum^m \Sb j=0 \\ j \not =i \endSb z_j \Biggl(\sum_{u''+|{\bold
v''}|=P} q_{(j,u'',{\bold v''})} a^{u''}{\bold b}^{\bold v''} \Biggr) \qquad
\endmultline \tag 14
$$
An easy rearrangement of the terms in $(14)$ yields the following
conclusion
$$
a^u {\bold b}^{\bold v}
\Biggl(z_i - \sum^m \Sb j=0 \\ j\not= i \endSb q_{(j,u, {\bold v})}
z_j \Biggr)
\in \biggl(a^l {\bold b}^{\bold w} \ :\ l+ |{\bold w}|=P,\
l \not = u \biggr),
$$
which, {\it a fortiori}, implies, by Remark~2.7, that
$$ 
z_i - \sum^m \Sb j=0 \\ j\not= i \endSb q_{(j,u, {\bold v})} z_j
\in \overline{(a, b_1, \ldots, b_{k-1})} \cap c(g). \tag 15
$$
However by our choice of $n$ and $t$ we have the following inclusions
$$
\overline{(a, b_1, \ldots, b_{k-1})} \cap c(g) \subseteq \overline{\m^n}
\cap c(g) \subseteq \m^t \cap c(g) \subseteq \m c(g),
$$
so that $(15)$ contradicts the minimality of the generators of $c(g)$.

If we pick $r_1, \ldots, r_{k-1}$ so that
$$
\left( \sum_{i=1}^{k-1} \frac{1}{r_i+1} \right)^{-1} > N,
$$
where  $N$ is the degree of $g$,
then the integer
$$
m \prod_{i=1}^{k-1} (r_i+1) -
\sum_{i=1}^{k-1} (N+r_i+1) (r_1+1) \cdots \widehat{(r_i+1)} \cdots (r_{k-1}+1)
$$
is strictly greater than $ \prod_{i=1}^{k-1} (r_i+1)$.

To complete the proof of Theorem~4.2, we use that $R$ is
approximately Gorenstein; hence there exists an
irreducible $\m$-primary ideal $Q$ of
$R$ such that the image of the ideal $c(g) \prod_{i=1}^{k-1} (a, b_i)^{r_i}$ in
$R/Q$ requires $m \prod_{i=1}^{k-1} (r_i+1)$ generators.

Let $g^*$, respectively $f_i^*$, denote the image of $g$,
respectively of the polynomial
$$
f_i = a^{r_i}t^{r_i}+a^{r_i-1}b_it^{r_i-1}+\cdots +ab_i^{r_i-1}t+b^{r_i},
$$
in $R/Q$. Then $c(f_i^*)$ is the image of $(a, b_i)^{r_i}$ in $R/Q$ and
the inequality in $(13)$ is satisfied with the $n$ of
Theorem~3.3 equal to $k-1$.
Hence there exists a polynomial $f_k^*$ in $(R/Q)[t]$
$($or polynomial $h$ in the notation of Theorem~3.3$)$ such that
$$
c(f_1^*) \cdots c(f_{k-1}^*)c(f_k^*) c(g^*) \not = 0
$$
while
$$
\left(\sum_{i=1}^{k-1} c(f_i^*g^*) c(f_1^*) \cdots \widehat{c(f_i^*)} \cdots c
(f_{k-1}^*)\right) c(f_k^*) + c(f_k^*g^*)c(f_1^*) \cdots c(f_{k-1}^*) = 0.
$$
This contradicts the fact that the polarized Dedekind--Mertens number of
$g^*$ is at most $k$ (see Remark~2.3). Hence $\mu(c(g)) = m \leq k$.

By Remark~1.2 and Theorem~2.5, we have
$\mu_R(c(g)) = \widetilde{\mu}_R(g) = \mu_R(g)$.
\QED
\enddemo

\medskip

\remark{Remark~4.3}
It would be interesting to know if Theorem~4.2
holds more generally without the hypothesis that the ring is
reduced. In the next section we present examples to show that
an assumption on the dimension of the ring is necessary in
Theorem~4.2.
\endremark

\head
5. One-dimensional examples
\endhead

In the presentation of the examples of this section  we use
Remark~2.2
that if $f$ and $g$ are polynomials in $R[t]$,
then the ideal $c(f)c(g)$ of $R$ is integral over $c(fg)$.

\example{Example~5.1} 
Let $F$ be a field and let $R$ be the Gorenstein subring $F[\![ s^3, s^4 ]\!]$
of the power series ring $\overline{R} = F[\![ s ]\!]$. Consider the polynomial
$$
g= s^7 + s^6t+ s^8 t^2 \in R[t].
$$
We claim that $\mu_R(g) = \widetilde{\mu}_R(g) = 2$ while $\mu_R(c(g)) = 3$.
\endexample
\demo{Proof}
We use the following  facts about the ring $R$.
\roster
\item"(1)"
The integral closure $\overline{I}$ of an ideal $I$ of $R$ is
$I\overline{R} \cap R$.

\item"(2)"
If $I$ is a nonzero ideal of $R$, then
$I\overline{R}  = s^n\overline{R}$ for
some nonnegative integer $n$. We have
$\mu_R(I) \le 3$ and $\mu_R(I)=3$ if and only
if $I$ also contains power series in $s$ of order $n+1$ and
$n+2$. In this case,
$I$ is integrally closed and
$I = I\overline{R} = s^n\overline{R}$,
i.e., $I$ is also
an ideal of $\overline{R}$, and $n \ge 6$.

\item"(3)"
If $\mu_R(I) = 3$ and $J$ is a non-zero ideal of $R$,
then $\mu_R(IJ)=3$, so $IJ$ is integrally closed.
\endroster

Statement (1) follows for example from \cite{ZS, Theorem~1, page~350}.
Since every nonzero ideal of $\overline{R}$ is of the form
$s^n\overline{R}$, the first sentence of Statement (2) is clear.
Since $R$ is a free module of rank 3 over the principal ideal
domain $F[\![s^3]\!]$, every nonzero ideal $I$ of $R$ is a free
$F[\![s^3]\!]$-module of rank 3. Therefore $\mu_R(I) \le 3$.
Suppose $\mu_R(I) = 3$ and $I = (h_1,h_2,h_3)R$.  If
$I\overline{R} = s^n\overline{R}$, then at least one of the
$h_i$ has order $n$ as a power series in $s$. We
may assume $h_1$ has order $n$. There exist $a,b \in F \subset R$
such that $h_2' = h_2 -ah_1$ and $h_3' = h_3 -bh_1$ have order greater than $n$.
We may assume that $h_2'$ has order less than or equal to that
of $h_3'$, and by subtracting from $h_3'$ a scalar multiple of $h_2'$,
we obtain $I = (h_1, h_2', h_3'')R$, where
the order of $h_2'$ is greater than $n$ and the order of $h_3''$ is
greater than the order of $h_2'$.
Moreover, if the order of $h_2'$ or $h_3''$ is $n+3$ or $n+4$,
we can subtract a scalar multiple of $s^3h_1$ or $s^4h_1$ to get
new generators of higher order. Thus we may assume $h_2'$ and
$h_3''$ have order different from $n+3$ and $n+4$.

Since $R$ contains all power series in $\overline{R}$ of order
greater than or equal to 6, the ideal $h_1R$ contains all
power series of order at least $n+6$.
If the order of $h_2'$ and $h_3''$ are
not $n+1$ and $n+2$, then the order of
$h_3''$ is at least $n+5$. However, if the order
of $h_3''$ is $n+5$, then
by subtracting from $h_3''$ a suitable
scalar multiple of $s^3h_2'$ or $s^4h_2'$,
we obtain a new minimal generator of
order at least $n+6$. This contradicts the fact that $\mu_R(I) = 3$.
Therefore the order of $h_2'$ and $h_3''$ must be $n+1$ and $n+2$.
Since $R$ contains no power series of order 1, 2 or 5, if
$\mu_R(I) = 3$, then $I$ has order at least 6.

Conversely, if $I$ has order $n$ and contains power series of
order $n+1$ and $n+2$, then $I$ is minimally generated by any
elements in $I$ with these orders, and in view of the fact that
$I$ contains all power series of order at least $n+6$, we see
that $I$ contains all power series of order greater than or
equal to $n$ and hence $I = s^n\overline{R}$.

Statement (3) follows from the characterization of 3-generated
ideals given in part (2).

\medskip

To establish Example~5.1,
we first show that $\mu_R(g)=2$. Let $f$ be an arbitrary polynomial in $R[t]$.
If $c(f)$ requires  three generators then the
equality $c(f)c(fg)=c(f)^2c(g)$ holds, as $c(f)c(fg)$ is
integrally closed  by (2) and $c(f)^2c(g)$ is integral over
$c(f)c(fg)$. On  the other hand, if $c(f)$ is principal then
the equality $c(f)c(fg) = c(f)^2c(g)$ also holds, as $f$ is a
Gaussian polynomial in the terminology of \cite{GV, HH1}.

In order to complete the proof it remains to consider the case in which
$\mu_R(c(f)) = 2$. We show
also in this case that the product $c(f)c(fg)$ requires three generators
and  hence is integrally closed.
Suppose first that $c(f)$ has order $n$ and contains a power series of
order $n+1$. Since $c(f)c(g)$ is integral over $c(fg)$, the ideal
$c(fg)$ contains a power series of order $n + 6$.  Suppose
$f = \sum a_it^i$ and let $j$ be minimal such that
$a_j$ as a power series in $s$ has order $n$ or $n+1$.
If $a_j$ has order $n$, then $c(fg)$ contains a power series of
order $n+7$. On the other hand, if $a_j$ has order $n+1$, then
$c(fg)$ contains a power series of order $n+8$.
In either case, since $c(f)$ contains power series of order $n$ and $n+1$
it follows that $c(f)c(fg)$ requires three generators and hence is
integrally closed.
The remaining possibility is that $c(f)$ has order $n$, contains a
power series of order $n+2$, and does not contain a power
series of order $n+1$.  In this case\footnote{In this case the
arrangement of the coefficients of $g$ does not matter.},
$c(fg)$ contains power series of order $n+6$ and $n+7$.
It follows that $c(f)c(fg)$ contains power series
of order $2n+6, 2n+7$ and $2n+8$ and therefore is integrally closed.

We show that the polarized Dedekind--Mertens number
$\widetilde{\mu}_R(g)=2$ as well. Let $f$ and $h$ be
two arbitrary polynomials in $R[t]$. If $c(f)$ is
principal, then $c(fg) = c(f)c(g)$, so
$c(fg)c(h) = c(f)c(g)c(h)$. Similarly, if $c(h)$ is principal,
we have $c(hg)c(f) = c(h)c(g)c(f)$.
On the other hand, if $\mu_R(c(f)) = 3$, then $c(hg)c(f)$ is integrally
closed and hence equal to $c(h)c(g)c(f)$. Similarly, if
$\mu_R(c(h)) = 3$, we have $c(fg)c(h) = c(f)c(g)c(h)$. It
remains to consider the case where $\mu_R(c(f)) = 2 = \mu_R(c(h))$.
In this case both $c(fg)c(h)$ and $c(hg)c(f)$ are integrally
closed, the argument being the same as that given in the
paragraph above.
\QED
\enddemo

\remark{Remark~5.2} 
We remark that the order of the coefficients of $g$ in Example~5.1 is
important. The polynomial
$g' = s^6+s^7t+s^8t^2$ has $c(g') = c(g)$, but $\mu_R(g') = 3$.
To see this, consider  the polynomial $f = s^6-s^7t$. We have
$fg' = s^{12}-s^{15}t^3$ so that $c(fg') = s^{12}R$. A direct calculation
shows that
$$
c(fg')c(f) = (s^{18}, s^{19})R \subsetneq (s^{18}, s^{19}, s^{20})R =
c(f)^2c(g').
$$
\endremark

\example{Setup~5.3} 
To generalize Example~5.1, let
$(R,\m)$ be a one-dimensional local Noetherian domain such that the integral
closure $(\overline{R},\overline{\m})$ of $R$ is
again local and is a finitely generated $R$-module.
Assume that the canonical injection $R/\m \hookrightarrow
\overline{R}/\overline{\m}$ is also surjective,
and let $s \in \overline{\m}$ be a generator
for the maximal ideal of the DVR $\overline{R}$.
Then $\m \overline{R} = s^e\overline{R}$,
where $e$ is the multiplicity of $R$. Since we are
assuming $\overline{R}$ to be a finitely generated $R$-module, the conductor of
$\overline{R}$ to $R$ is a nonzero ideal of
$\overline{R}$, so it is of the form $s^c\overline{R}$ for some
nonnegative integer $c$.

If $I$ is a nonzero ideal of $R$, then
$I\overline{R} = s^n\overline{R}$ for some nonnegative
integer $n$. To better measure and compare ideals of $R$, we associate
with $I$ a subset $\gamma(I)$ of the positive integers less than $e$,
defined as follows:
$$
\gamma(I) = \{ i \in {\Bbb N} \,|\, 1 \leq i < e \text{ and }
\exists a \in I \text{ with } a\overline{R} = s^{n+i}\overline{R}\}.
$$

In analogy with the observations made at the beginning of the
proof of Example~5.1, we have the following facts about $R$.
\roster
\item"(1)"
The integral closure $\overline{I}$ of an ideal $I$ of $R$ is
$I\overline{R} \cap R$.
A nonzero ideal $I$ of $R$ has $\gamma(I) = \{ 1, 2, \ldots, e-1 \}$ if and
only if $I = I\overline{R}$; in this case, $I$ is
integrally closed and is contained in
the conductor of $\overline{R}$ to $R$.

\item"(2)"
For every ideal $I$ of $R$ we have $\mu_R(I) \le e$ and $\mu_R(I) \geq 1 +
|\gamma(I)|$, where $|\gamma(I)|$ denotes the cardinality of $\gamma(I)$.
A nonzero integrally closed ideal $I$ of $R$ contained in $s^c\overline{R}$
has $\gamma(I) = \{ 1, 2, \ldots, e-1 \}$.
In particular, every nonzero integrally closed ideal $I$ of $R$ that is
contained in the conductor of $\overline{R}$ to $R$
has\footnote{In this more general
setting, as contrasted with Example~5.1, there exist rings $R$ with
multiplicity $e$ and non-integrally closed ideals $I$ of $R$ with
$\mu_R(I)=e$. For example, consider the subring $R=F[\![s^3, s^5]\!]$ of the
formal power series ring $\overline{R}=F[\![ s ]\!]$. Then $I=(s^3, s^5)^3 =
(s^9, s^{11}, s^{13})$ has $\mu_R(I)=3$,
but $I$ is not integrally closed.
Note that $\gamma(I)=\{2\}$.}
precisely $\mu_R(I) = e$.

\item"(3)"
If $I$ and $J$ are nonzero ideals of $R$ then $\gamma(IJ)$ contains
$\gamma(I)$ and $\gamma(J)$ as well as all $i+j \leq e-1$ where
$i \in \gamma(I)$ and $j \in \gamma(J)$.
In particular, if $|\gamma(I)|=e-1$ and $J$ is a nonzero ideal of $R$
then $|\gamma(IJ)|=e-1$, so $IJ$ is integrally closed.

\item"(4)"
Let $f, g \in R[t]$ be nonzero polynomials and let $m$ be a positive integer
less than $e$. If $\{ 1, \ldots, m \} \cap \gamma(c(f)) =
\varnothing$ and $\{ 1, \ldots, m \} \subseteq \gamma(c(g))$, then
$\{ 1, \ldots, m \} \subseteq \gamma(c(fg))$.
\endroster

\noindent
In view of (2), Theorem~2.5 implies that
$\widetilde{\mu}_R(g) \le e$ for every polynomial $g \in R[t]$.
\endexample

\example{General Example 5.4} 
With notation as in Setup~5.3, assume that $R$ has multiplicity
$e(R) = e \ge 3$, and let $s^c\overline{R}$ be the conductor
of $\overline{R}$ to $R$.
Specific examples for $R$ are, for instance, the subrings of the formal power
series ring $\overline{R}=F[\![ s ]\!]$ of the form: $R=F[\![ s^e,
s^{e+1} ]\!]$ or $R=F[\![ s^e, s^{e+1}, \ldots, s^{2e-1} ]\!]$. In the first
case $R$ is a complete intersection (therefore Gorenstein) and $c = (e-1)e$,
while in the second case $c = e$ and $R$ fails to be Gorenstein.

In analogy with Remark~5.2, consider the polynomial
$$
g'= s^c + s^{c+1}t + s^{c+2}t^2 + \ldots + s^{c+e-1}t^{e-1} \in R[t].
$$
Let $f=s^c - s^{c+1}t$. Then $fg'= s^{2c}-s^{2c+e}t^e$, and hence
$c(fg')=s^{2c}R$. By Setup~5.3(1), we have $c(g')=s^c\overline{R}$
and by Setup~5.3(3) we have that $c(g')c(f)^k$ is integrally closed for each
positive integer $k$. Therefore the smallest positive integer $k$ for which
$$
c(fg') c(f)^{k-1} = c(f)^k c(g')
$$
is the smallest integer $k$ such that $c(fg') c(f)^{k-1}$ is integrally closed.
Since $c(fg') = s^{2c}R$ is principal in $R$, $c(fg') c(f)^{k-1}$ is integrally
closed precisely when $c(f)^{k-1}$ is integrally closed. As $c(f)^{k-1}$
is generated by $k$ elements, it is, in view of Setup~5.3(2),
integrally closed only for $k \geq e$. Therefore $\mu_R(g') = e =
\mu_R(c(g'))$.

\smallskip

On the other hand, the polynomial
$$
g = s^{c+1} + s^ct + s^{c+2}t^2 + \ldots + s^{c+e-1}t^{e-1}
$$
has $c(g) = c(g')$, but $\mu_R(g) \leq e-1 < e=\mu_R(c(g))$. To justify
this assertion, let $f \in R[t]$ be a nonzero polynomial; the following
two cases are possible: (a) $1 \not\in \gamma(c(f))$; (b) $1 \in \gamma(c(f))$.

\roster
\item"(a)"
Choose $1 \leq m \leq e-1$ maximal with the
property that $\{ 1, \ldots, m\} \cap \gamma(c(f)) = \varnothing$.
As $\gamma(c(g))= \{ 1, \ldots, e-1 \}$, by Setup~5.3(4)  we
have $\{1, \ldots, m \} \subseteq \gamma(c(fg))$.
By Setup~5.3(3), $\gamma(c(fg)c(f)^{k-1})$ contains each of the
integers $1, \ldots, km+k-1$ that is less than $e$.
Hence, for $k$ a positive integer such
that $km+k-1 \geq e-1$, the ideal $c(fg)c(f)^{k-1}$ is integrally closed
and therefore equal to $c(f)^kc(g)$.
Choosing $k = \lceil \frac{e}{m+1} \rceil$ guarantees that
$c(fg)c(f)^{k-1}=c(f)^kc(g)$. Finally, notice that
$$
k=\left\lceil \frac{e}{m+1} \right\rceil < e.
$$
If not, $\lceil \frac{e}{m+1} \rceil \geq e$ or, equivalently,
$\frac{e}{m+1} > e-1$. But this last inequality yields the contradicting
conclusion $1 > m(e-1)$.

\item"(b)"
If $\gamma(c(f)) = \{1, \ldots, e-1 \}$, then by Setup~5.3(4),
$c(fg)c(f)$ is integrally closed. Hence $c(fg)c(f) = c(f)^2c(g)$ in
this case. In the case that remains, there exists a positive
integer $m$ less than $e-1$ such that
$ \{ 1, \ldots, m \} \subseteq \gamma(c(f))$  while
$m+1 \not\in \gamma(c(f))$.

The following fact that depends on the specific
arrangement of the coefficients of $g$ implies that
$\{ 1, \ldots, m+1 \} \subseteq \gamma(c(fg)c(f))$.

\medskip

\itemitem{(b.1)}
Let $f \in R[t]$ be a nonzero polynomial and let $m$ be a positive integer
less than $e-1$. If $ \{ 1, \ldots, m \} \subseteq \gamma(c(f))$ and
$m+1 \not\in \gamma(c(f))$, then $\{ 1, \ldots, m+1 \} \cap \gamma(c(fg))
\not= \varnothing$.

\medskip

By Setup~5.3(3), for an integer $k \ge 2$,
$\gamma(c(fg)c(g)^{k-1})$ contains each of the integers $1, \ldots, (k-1)m+1$
that is less than $e$.
Hence, for $k$ a positive integer such
that $(k-1)m+1 \geq e-1$, the ideal $c(fg)c(f)^{k-1}$ is integrally closed
and therefore equal to $c(f)^kc(g)$.
Choosing $k = \lceil \frac{e+m-2}{m} \rceil$ guarantees that
$c(fg)c(f)^{k-1}=c(f)^kc(g)$. Finally, notice that
$$
k=\left\lceil \frac{e+m-2}{m} \right\rceil < e.
$$
If not, $\lceil \frac{e+m-2}{m} \rceil \geq e$ or, equivalently,
$e+m-2 > (e-1)m$. But this last inequality yields the contradicting
conclusion $2(m-1) > e(m-1)$.
\endroster

\noindent
Putting (a) and (b) together we conclude that $\mu_R(g) \leq e-1$, as
claimed.

\smallskip

To show the polarized Dedekind--Mertens number $\widetilde{\mu}_R(g)$
is also at most $e-1$, we need to show
for polynomials $f_1, \ldots, f_{e-1} \in R[t]$ that we have
$$ 
\sum_{i=1}^{e-1} c(f_ig) c(f_1) \cdots \widehat{c(f_i)} \cdots c(f_{e-1})
=
c(f_1) \cdots c(f_{e-1}) c(g). \tag 16
$$
Since
$c(f_1) \cdots c(f_{e-1}) c(g)$ is integral over
$c(f_ig) c(f_1) \cdots \widehat{c(f_i)} \cdots c(f_{e-1})$ for each $i$, it
suffices to show one of these ideals is integrally closed in order
to establish $(16)$.
To do this, we consider the following cases.

\roster
\item"(1)"
Suppose  $ \gamma(c(f_i)) = \{1,\ldots, e-1 \}$ for some $i$.
Then Setup~5.3(3) implies that the ideal
$c(f_jg) c(f_1) \cdots \widehat{c(f_j)} \cdots c(f_{e-1})$ is
integrally closed for every $j \ne i$, so $(16)$ holds
in this case. Hence we assume
$|\gamma(c(f_i))| \le e-2$ for each $i$. Applying
either Setup~5.3(4) or (b.1), it follows
that $\gamma(c(f_ig)) \ne \varnothing$ for each $i$.

\item"(2)"
Suppose  $\gamma(c(f_i)) = \varnothing$ for some $i$. Then
Setup~5.3(4) implies that $|\gamma(c(f_ig))| =  e-1 $.
By Setup~5.3(3),
$c(f_ig) c(f_1) \cdots \widehat{c(f_i)} \cdots c(f_{e-1})$ is
integrally closed,  so $(16)$ holds in this case. Hence we assume
each $\gamma(c(f_i)) \ne \varnothing$.

\item"(3)"
Suppose $1 \in \gamma(c(f_i))$ for every $i$, $1 \le i \le e-1$.
In view of (1), Setup~5.3 implies that the ideal
$c(f_1g)c(f_2) \cdots c(f_{e-1}) = I$ has $|\gamma(I)| = e-1$ and
hence is integrally closed. Therefore $(16)$ holds in this
case. Hence we assume $1 \not\in \gamma(c(f_i))$ for some $i$.

\item"(4)"
Let  $m$ be the  positive integer maximal with the property
that for every $i$, $1 \le i \le e-1$,
there exists a positive integer $k_i \in \gamma(c(f_i))$
with $k_i \le e-m$. In view of (2), there exists such an integer $m$
and in view of (3), $m \le e-2$.
By the maximality of $m$, we have
$\gamma(c(f_i)) \subseteq  \{e-m, e-m+1,  \cdots ,  e-1 \}$ for
some $i$.
Setup~5.3(4) implies that $\gamma(c(f_ig))$ contains
$\{1, \ldots, e-m-1 \}$. Since for each $j \ne i$,
there exists a positive integer
$k_j \in \gamma(c(f_j))$ with $k_j \le e-m$,
it follows from Setup~5.3 that  the ideal
$I = c(f_ig) c(f_1) \cdots \widehat{c(f_i)} \cdots c(f_{e-1})$
has $|\gamma(I)| = e-1$ and hence is integrally closed.
Therefore $(16)$ holds in general.

We conclude that $\widetilde{\mu}_R(g) \le e-1$.
\endroster
\endexample

\remark{Remark~5.5} \noindent \roster
\item"(1)"
For the polynomial $g$ in General Example~5.4 and $f=s^c-s^{c+1}t$,
a simple computation shows that
$c(fg)=(s^{2c}, s^{2c+1}, s^{2c+2})R$. Since $\gamma(c(f))=\{1\}$ and
$\gamma(c(fg))=\{1, 2\}$, for $k$ a positive integer
the ideal $c(fg)c(f)^{k-1}$ is integrally closed (and hence equal to
$c(f)^kc(g)$) if and only if $k \ge e-2$. Hence
the Dedekind--Mertens number
$\mu_R(g)$,  as well as the polarized Dedekind--Mertens
number $\widetilde{\mu}_R(g)$, is either $e-1$ or $e-2$.

\item"(2)"
If $e =3$, then $\mu_R(g) = e-1 = 2 = \widetilde{\mu}_R(g)$;
for otherwise $g$ would be Gaussian, and it is known
\cite{HH1, Theorem~1.5}
that over a Noetherian
local domain a Gaussian polynomial has principal content ideal.

\item"(3)"
To obtain an example where  $\mu_R(g) < \mu_R(c(g)) - 1$,
we modify the $g$ of General Example~5.4 as follows:
With $R$ as in General Example~5.4 and $e = 5$, let
$$
g = s^{c+1} + s^ct + s^{c+2}t2 + s^{c+4}t^3 + s^{c+3}t^4.
$$
Then $\mu_R(g) \le 3$, while $\mu_R(c(g)) = 5$.
To see that $\mu_R(g) \le 3$, we examine cases similar to
what is done above. The `new' case is where
$\gamma(c(f)) = \{ 1 \}$. In this case with our modified
$g$, the set  $\gamma(c(fg))$ contains either
1 or 2 and also contains either 3 or 4.
\endroster
\endremark

\remark{Remark~5.6}
If $R$ is a Noetherian domain of dimension at
least $2$,  it would be interesting to know
whether, in situations where  the hypotheses of Theorem~4.2
are not satisfied,  there
exist polynomials $g$ over $R$ with $\mu(g) < \mu(c(g))$.
For example, if $R$ is the polynomial ring $k[x, y]$ with $k$ a
field and $g \in R[t]$ is a polynomial such that $(c(g))= (x, y)^3$,
could  it happen that $\mu(g) \leq 3$?
\endremark

\enddocument